\documentclass[12pt,oneside,reqno]{amsart}
\newtheorem{theorem}{Theorem}

\newtheorem{definition}[theorem]{Definition}

\usepackage{xcolor}

\numberwithin{equation}{section}
\title[Pursuit Game Described by Infinite System of DE]{Pursuit Game for an Infinite
System of First-Order Differential Equations with Negative Coefficients}

\author[G.I. Ibragimov]{Gafurjan Ibragimov}
\address[Gafurjan Ibragimov]{Department of Mathematics and Institute for Mathematical
Research, Universiti Putra Malaysia, Serdang, Malaysia} \email{ibragimov@upm.edu.my}
\author[M. Ferrara]{Massimiliano Ferrara}
\address[Massimiliano Ferrara]{Department of Law, Economics and Human Sciences and Decisions Lab, University Mediterranea of
Reggio Calabria, Italy \& ICRIOS - The Invernizzi Centre for Research in Innovation, 
Organization, Strategy and Entrepreneurship Bocconi University - Department of Management and Technology, Italy} 
\email{massimiliano.ferrara@unirc.it, massimiliano.ferrara@unibocconi.it}
\author[I.A. Alias]{Idham Arif Alias}
\address[Idham Arif Alias]{Department of Mathematics and Institute for Mathematical Research,
Universiti Putra Malaysia, Serdang, Malaysia} \email{idham\_2@upm.edu.my}
\author[Mehdi Salimi]{Mehdi Salimi}
\address[Mehdi Salimi]{Department of Law, Economics and Human Sciences and Decisions Lab, University Mediterranea of
Reggio Calabria, Italy \&
Center for Dynamics and Institute for Analysis, Department
of Mathematics, Technische Universit{\"a}t Dresden, Germany}
\email{mehdi.salimi@unirc.it, mehdi.salimi@tu-dresden.de}

\begin{document}
\begin{abstract}
In this paper we study a linear pursuit differential game described by an infinite system of first-order differential equations in Hilbert space. The control functions of players are subject to geometric constraints. The pursuer attempts to bring the system from a given initial state to the origin for a finite time and the evader's purpose is opposite. We obtain a guaranteed pursuit time and construct a strategy for pursuer.
\end{abstract}
\maketitle

\section{Introduction}
A game has players, strategies and scores to explain why the players win or lose.
Differential games are games which are modeled with differential equations and were initiated by Isaacs \cite{isa}. Pursuit--evasion games are common examples of games that are abstract models for pursuers who try to catch evaders who
are running away. Differential games and pursuit--evasion problems
are investigated by many authors. For example see, Blagodatskih and Petrov \cite{bla}, Krasovskii \cite{kra},  Pashkov and Terekhov
\cite{pas}, Petrov \cite{pet1,petrov}, Petrosyan \cite{pet},
Pontryagin \cite{pon}, Pshenichii \cite{psh}, Rikhsiev \cite{rik} and Rzymowski
\cite{rzy}. Finding the value of the game and identifying optimal
strategies of players are two interesting subjects in the study of
pursuit evasion problems. Such issues for many players with
a variety of constraints have been studied in \cite{fer}, \cite{ibr0}, \cite{ibr1}, \cite{ibr2}, \cite{ibr3}, \cite{ibr4}, \cite{salsem}, \cite{salfer} and \cite{saldyn}.

In \cite{ibr1} Ibragimov investigates a pursuit--evasion game of optimal
approach of countably many pursuers to one evader in a Hilbert
space with geometric constraints on the controls of the players.
Ibragimov and Salimi \cite{ibr3} workedon a similiar game
for inertial players with integral constraints with the
assumption that the control resource of the evader is less than
that of each pursuer. Evasion of evader from many pursuers in simple motion
differential games was developed by
Ibragimov et al.\ in \cite{ibr4} as well. This class of
differential games occurs naturally as a reformulation for solving
controlled distributed systems described by a parabolic equation
as in \cite{sat} (see also \cite{cher,ibr2} for more details). In \cite{satt} Satimov and Tukhtasinov study pursuit and evasion problems for controlled equations of the parabolic type. The control
parameters occur on the right-hand side of the equations as additive terms. They consider various
cases of control constraints. For some cases, they single out pairs of sets of initial states such that
capture is guaranteed if the initial point belongs to the first set and evasion of the terminal set is
guaranteed if the initial point belongs to the second set.

In \cite{ibbtai} Ibragimov et al. study a differential game of optimal approach of finite or countable
number of pursuers with one evader in the Hilbert space $l_2$. In their research, they found a formula for the value of the
game and constructed explicitly optimal strategies of the players. Important point to
note was that the energy resource of any pursuer needs not be greater than that of
the evader.

In the work of  Huang et al. \cite{Hao}, a game of multiple pursuers cooperating to capture a single evader in a bounded, convex polygon was studied in the plane. The main result of that paper is construction of a decentralized, guaranteed pursuit strategy where the pursuers cooperatively minimize the area of the evaders Voronoi cell by independently controlling each pursuers shared Voronoi boundary with the evader. In the differential game of many pursuers in a planar domain studied by Zhengyuan et al. \cite{Zhe} pursuit strategies are constructed based on Voronoi partition as well.

In the paper of Salimi and Ferrara \cite{salfer}, authors consider a finite time pursuit--evasion game in which a finite or
countable number of pursuers pursue a single evader. The control
functions of players are subjected to integral constraints.They introduce the value of the
game and identify optimal strategies of the pursuers. Azamov and Ruziboev \cite{Azru} considered the time-optimal problem for a controlled system with evolution-type distributed parameters. They obtained an upper estimate for the optimal transition time into the zero state.

In order to motivate the study of countably many pursuers in
contrast to the classical problem of finitely many pursuers, we
would like to draw a comparison with optimal control theory
\cite{pesch}. In order to design a controller for an
infinite-dimensional problem, e.g.\ a partial differential
equation, one often follows the `approximate-then-design' method,
i.e.\ first the infinite-dimensional problem is approximated,
e.g.\ by a Galerkin approximation, and then a controller is
designed for the resulting finite-dimensional approximation.
Finally the finite-dimensional controller or strategy is shown to
work also for the original problem \cite{elfarra}. A dual approach
which is still in its infancy is the `design-then-approximate'
method, i.e.\ first an infinite-dimensional controller or strategy
for the infinite-dimensional problem is designed, and then
approximated by a finite-dimensional controller which is then
proved to work also for the original problem \cite{burns}. The
results in this paper are in the spirit of the
`design-then-approximate' method and in case one is interested in
finite-dimensional approximations or implementations of our
methods, further investigations are required.

\section{Statement of problem}\label{sec:2}
Consider the Hilbert space
\[
l^2=\left\{\alpha=(\alpha_1,\alpha_2,...)|\ \sum\limits^\infty_{k=1} \alpha^2_k < \infty \right\}
\]
with inner product and norm defined by
\[
\langle\alpha, \beta\rangle = \sum\limits_{k = 1}^{\infty}|\alpha_k\beta_k|, \ \alpha, \ \beta \in l^2, \ \ ||\alpha|| =
\left(\sum\limits_{k = 1}^{\infty} \alpha_k^2 \right)^{1/2}
\]
and the following differential game described by the equations
\begin{equation}\label{a}
\dot{z_i}=-\lambda_iz_i+u_i-v_i, \quad i=1,2,\dots; \ \ z_i(0)=z_{i0},
\end{equation}
where $z_i, u_i, v_i \in \mathbb{R}$, $u=(u_1,u_2, \dots)$ and $v=(v_1,v_2, \dots)$ are control parameters of pursuer and evader respectively, and $\lambda_k$, $k=1,2,...$, are positive numbers. It is assumed that $z_0=(z_{10},z_{20}, \dots) \ne 0$.

\begin{definition}
Let $T$ be an arbitrary number. A vector function $u(t)=(u_1(t),u_2(t),\dots)$, $||u(t)||\leq\rho, 0\leq t\leq T,$ is called admissible control of the pursuer, where $\rho$ is a given positive number
\end{definition}

\begin{definition}
A vector function $v(t)=(v_1(t),v_2(t),\dots)$, $||v(t)||<\sigma$, $0\leq t\leq T$ is called admissible control of the evader, where $\sigma$ is a given positive number
\end{definition}

\begin{definition}
A function of the form $u(v)=(u_1(v_1),u_2(v_2),\dots)$ with measurable coordinates of $u_i(\cdot)\in\mathbb{R}$, that satisfies the condition $||u(v)||\leq \rho$ for any admissible control of evader $v=v(t),0\leq t\leq T$, is called strategy of pursuer
\end{definition}

Pursuit is started from the initial positions $z_0=(z_{10},z_{20},\dots )$ at time $t=0$ where \par\noindent$z_{i0}\in \mathbb{R}$,$i=1,2,\dots $.\par 

If we replace the parameters $u_i,v_i$ in the equation (\ref{a}) by some measurable functions $u_i(t), v_i(t), 0\leq t\leq T$, then it follows from the theory of differential equations that the initial value problem (\ref{a}), $z_0=(z_{10},z_{20},\dots )$ has a unique solution on the time interval $[0,T]$.\par 

The solution $$z(t)=(z_1(t),z_2(t),\dots),\qquad 0\leq t\leq T$$
of infinite system of differential equations (\ref{a}) is considered in the space of functions $f(t)=(f_1(t),f_2(t),\dots )\in l^2$ with absolutely continuous coordinates $f_i(t)$ defined on the interval $0\leq t\leq T$.

\begin{definition}
A number $T_0, T_0\leq T$, is called a guaranteed pursuit time if there exists a strategy of pursuer such that for any control of the evader, the solution of the initial value problem (\ref{a}) , $z_0=(z_{10},z_{20},\dots )$, $z(t)=(z_1(t),z_2(t),\dots)$ equals zero, at some $T$, $0\leq T\leq T_0$, i.e. $z_i(T)=0$ for all $i=1,2,\dots $
\end{definition}

\noindent{\bf{Problem} :} Find a guaranteed pursuit time in the game (\ref{a})-($z_0=(z_{10},z_{20},\dots) $).
\section{Main result}

We assume that $\rho > \sigma$ and define the strategy for the pursuer as following:
\begin{align}\label{b}
&u_i(t) = \left\{ \begin{array}{ll}
                -\frac{z_{i0}}{||z_0||}(\rho-\sigma)+v_i(t), i=1,2,..., & 0\leq t\leq T_i, \\
                v_i(t), & t>T_i \\
                \end{array} \right.&&
\end{align}
where $$T_i=\frac{1}{\lambda_i}\ln(\frac{\lambda_i||z_0||}{\rho-\sigma}+1).$$
The above strategy is admissible, indeed using the Minkowskii inequality we get,
\begin{align*}
||u||&=\Bigg(\sum_{i=1}^\infty \bigg(-\frac{z_{i0}}{||z_0||}(\rho-\sigma)+v_i\bigg)^2\Bigg)^\frac{1}{2}\\
&\leq \bigg(\sum_{i=1}^\infty \frac{z_{i0}^2(\rho-\sigma)^2}{||z_0||^2}\bigg)^\frac{1}{2}+\bigg(\sum_{i=1}^\infty v_i^2\bigg)^\frac{1}{2}\\
&\leq \rho-\sigma+\sigma=\rho.
\end{align*}

\noindent Thus, $||u||\leq \rho$ and the strategy is admissible. \par \medskip
\noindent In this part we are going to prove that,
\begin{align}\label{c}
z_i(t)=0 \quad\mbox{for all}\quad t\geq T_i.
\end{align}
\noindent Indeed, by (\ref{b})  we have
\begin{align*}
z_i(t)&=e^{-\lambda_it}\bigg[z_{i0}+\int_0^te^{\lambda_is}(u_i(s)-v_i(s))ds\bigg]\\
&=e^{-\lambda_it}\bigg[z_{i0}-\int_0^{T_i}e^{\lambda_is}\frac{z_{i0}}{||z_0||}(\rho-\sigma)ds\bigg]=0.
\end{align*}

Since
\begin{equation*}
\int_0^{T_i} e^{\lambda_is}ds=\frac{e^{\lambda_iT_i}-1}{\lambda_i}=\frac{||z_0||}{\rho-\sigma}.
\end{equation*}

Thus, by (\ref{c}) the number
\begin{equation*}
T=\underset{i}{\sup}\ T_i=\underset{i}{\sup}\ \frac{1}{\lambda_i}\ln\big(\frac{\lambda_i||z_0||}{\rho-\sigma}+1\big)
\end{equation*}
is a guaranteed pursuit time.\par \medskip

\noindent Now we show that
\begin{equation}\label{d}
T=\frac{1}{\lambda}\ln\big(\frac{\lambda||z_0||}{\rho-\sigma}+1\big)
\end{equation}
where $\lambda=\underset{i}{\inf}\ \lambda_i$. To this end we consider the function
\begin{equation*}
f(x)=\frac{\ln\big(\frac{||z_0||}{\rho-\sigma}x+1\big)}{x},\quad x>0.
\end{equation*}

It is not difficult to see that
\begin{equation*}
f'(x)=\frac{1}{x^2}g(x),\quad g(x)=1-\frac{1}{\frac{||z_0||}{\rho-\sigma}x+1}-\ln\big(\frac{||z_0||}{\rho-\sigma}x+1\big),
\end{equation*}
and
\begin{equation*}
g'(x)=-\frac{\Big(\frac{||z_0||}{\rho-\sigma}\Big)^2x}{\Big(\frac{||z_0||}{\rho-\sigma}x+1\Big)^2}<0.
\end{equation*}

Hence $g(x)$ is decreasing. Since $g(0)=0$, therefore $g(x)<0,\ x>0$. Therefore, $f'(x)<0, \ x>0$. Consequently $f(x)$ is decreasing. Thus, $T$ is defined by formula (\ref{d}).\par \medskip

\section{Discussion and Conclusion}

We studied a linear pursuit differential game described by an infinite system of first-order differential equations in Hilbert space $l^2$. The control functions of players are subjected to geometric constraints. The pursuer attempts to bring the system from a given initial state to the origin for a finite time and the evader's purpose is opposite. We obtained a guaranteed pursuit time and constructed a strategy for pursuer.

In the differential game studied by Satimov and Tukhtasinov \cite{sat} the numbers $\lambda_1,\lambda_2,\dots,$ satisfy the inequality $0<\lambda_1\leq\lambda_2\leq \dots \rightarrow \infty $ and guaranteed pursuit time is $T_0=\frac{||z_0||}{\rho-\sigma}$. For that differential game guaranteed pursuit time $T$ defined by the formula (\ref{d}) is
\begin{equation*}
T=\frac{1}{\lambda_1}\ln\big(\frac{\lambda_1||z_0||}{\rho-\sigma}+1\big)<T_0=\frac{||z_0||}{\rho-\sigma}.
\end{equation*}
In summary, we stress the following two advantages of the present work
\begin{enumerate}
\item Guaranteed pursuit time $T$ defined by (4) improves that of the work of Satimov, Tukhtasinov \cite{sat}.
\item The numbers $\lambda_1,\lambda_2,\dots,$ needn't satisfy the relations $\lambda_1\leq\lambda_2\leq \dots$. They are assumed to be any positive numbers.
\end{enumerate}


\begin{thebibliography}{99}

\bibitem{Azru} 
Azamov A.A., Ruziboev M.B., The time-optimal problem for evolutionary partial
differential equations, Journal of Applied
Mathematics and Mechanics, 77(2): 220-224, (2013)

\bibitem{bla} Blagodatskih, A.I., Petrov, N.N., Conflict Interaction
Controlled Objects Groups, Izhevsk: Udmurt State University (in
Russian), (2009)

\bibitem{burns} Burns, J.A., King, B.B., A reduced basis approach to the design of low-order feedback
controllers for nonlinear continuous systems, Journal of Vibration
and Control, 4, 297-323, (1998)

\bibitem{cher} Chernous'ko, F.L., Bounded controls in distributed-parameter systems, Journal of Applied
Mathematics and Mechanics, 56(5), 707-723, (1992)

\bibitem{elfarra} El-Farra, N.H, Armaou, A., Christofides, P.D., Analysis and control of
 parabolic PDE systems with input constraints, Automatica, 39, 715-725,
 (2003)

\bibitem{fer} Ferrara, M., Ibraimov, G., Salimi, M., Pursuit-evasion game
of many players with coordinate-wise integral constraints on a
convex set in the plane, AAPP | Atti della Accademia Peloritana
dei Pericolanti Classe di Scienze Fisiche, Matematiche e Naturali,
95(2), A6, (2017)

\bibitem{Hao} Huang, H., Zhangy, W., Ding, J., Stipanovi\'c, D.M., Tomlin, C.J., Guaranteed Decentralized Pursuit-Evasion in the Plane with Multiple Pursuers, 50th IEEE Conference on Decision and Control and European Control Conference (CDC-ECC) Orlando, FL, USA, Dec. 12-15, (2011)

\bibitem{ibr0}Ibragimov, G.I., A game of optimal pursuit of one object by
several, Journal of Applied Mathematics and Mechanics, 62(2),
187-192, (1998)

\bibitem{ibr1} Ibragimov, G.I., Optimal pursuit with countably many pursuers
and one evader, Differential Equations, 41(5), 627-635, (2005)

\bibitem{ibr2}  Ibragimov G.I., The optimal pursuit problem reduced to an
infinite system of differential equations, Journal of Applied
Mathematics and Mechanics, 77, 470-476, (2013)

\bibitem{ibbtai} Ibragimov, G., Norshakila, A.R., Kuchkarov, A., Fudziah, I., Multi Pursuer Differential Game of Optimal Approach with Integral
Constraints on Controls of Players, Taiwanese Journal of
Mathematics 19(3), 963-976, (2015)

\bibitem{ibr3} Ibragimov, G.I., Salimi, M., Pursuit-evasion differential game with
 many inertial players, Mathematical Problems in Engineering, vol. 2009, Article ID
653723, 15 pages, (2009)

\bibitem{ibr4} Ibragimov, G.I., Salimi, M., Amini, M., Evasion from many pursuers
in simple motion differential game with integral constraints,
European Journal of Operational Research, 218, 505-511, (2012)

\bibitem{isa} Isaacs, R.,  Differential Games, John Wiley \& Sons, New York, NY, USA,
(1965)

\bibitem{kra} Krasovskii, N.N., Control of a Dynamical System, Nauka, Moscow,
(1985)

\bibitem{pas} Pashkov, A.G., Terekhov, S.D., On a game of optimal pursuit of
one object by two objects, Prikl. Mat. Mekh., 47(6), 898-903,
(1983)

\bibitem{pesch} Pesch, H.J., Solving optimal control and pursuit-evasion game problems of high complexity,
 Computational Optimal Control, 115, 43-61, (1994)

\bibitem{pet1} Petrov, N.N., Simple Group Pursuit Subject to Phase Constraints and Data Delay, Journal of Computer and Systems Sciences International, 57(1): 37-42, (2018)

\bibitem{petrov} Petrov, N.N., About one group pursuit problem with phase
constraints, Prikladnaya Matematika i mechanika, 6, 1060-1063,
(1988)

\bibitem{pet} Petrosyan, L.A., Differential Pursuit Games, Izdat. Leningrad. Univ., Leningrad,
(1977)

\bibitem{pon} Pontryagin, L.S., Selected Works,  Moscow: MAKS Press.,
(2004)

\bibitem{psh} Pshenichii, B.N., Simple pursuit by several objects, Cybernetics and Systems
Analysis, 12(3): 145-146, (1976)

\bibitem{rik} Rikhsiev, B.B., The Differential Games with Simple Motions, Tashkent:
Fan, (1989)

\bibitem{rzy} Rzymowski, W., Evasion along each trajectory in differential games with many pursuers,
 Journal of Differential Equations, 62(3), 334-356, (1986)

\bibitem{salsem} Salimi, M., A research contribution on an evasion problem,
SeMA Journal, 75(1), 139-144, (2018)

\bibitem{salfer} Salimi, M., Ferrara, M., Differential game of optimal pursuit of one evader by many pursuers,
International Journal of Game Theory, 48(2), 481-490, (2019)

\bibitem{saldyn} Salimi, M., Ibragimov, G., Siegmund, S., Sharifi, S., On a fixed
duration pursuit differential game with geometric and integral
constraints, Dynamic Games and Applications 6(3), 409-425, (2016)

\bibitem{sat} Satimov, N.Yu., Tukhtasinov, M., Game problems on a fixed
interval in controlled first-order evolution equations,
Mathematical notes, 80(3-4), 578-589, (2006)

\bibitem{satt} Satimov, N.Yu., Tukhtasinov, M., On Some Game Problems
for First-Order Controlled Evolution Equations,
Differential Equations, 41(8), 1169-1177, (2005)

\bibitem{Zhe} Zhou, Z., Zhang, W., Ding, J., Huang, H., Stipanovi\'c, D.M., Tomlin, C.J., Cooperative pursuit with Voronoi partitions, Automatica, 72, 64-72, (2016)

\end{thebibliography}
\end{document}